\documentclass[12pt]{article}
\usepackage{amsmath}
\usepackage{amssymb}
\usepackage{amsfonts}
\usepackage{eucal}
\usepackage[usenames]{color}
\usepackage{graphicx}
\numberwithin{equation}{section}
\newtheorem{Theorem}{Theorem}[section]
\newtheorem{Proposition}[Theorem]{Proposition}

\newtheorem{Remark}[Theorem]{Remark}

\begin{document}
\title{On the first-passage time of an integrated Gauss-Markov process}
\author{Mario Abundo\thanks{Dipartimento di Matematica, Universit\`a  ``Tor Vergata'', via della Ricerca Scientifica, I-00133 Rome, Italy.
E-mail: \tt{abundo@mat.uniroma2.it}}
}
\date{}
\maketitle

\begin{abstract}
\noindent It is considered the integrated process $X(t)= x + \int
_0^t Y(s) ds ,$ where $Y(t)$ is a Gauss-Markov process starting
from $y.$ The first-passage time (FPT) of $X$ through a constant
boundary and the first-exit time of $X$ from an interval $(a,b)$
are investigated, generalizing some results on FPT of integrated
Brownian motion. An essential role is played by a useful
representation of $X,$
which allows
to reduces the FPT of $X$ to that of a time-changed Brownian
motion. Some explicit examples are reported; when theoretical
calculation is not available, the quantities of interest are
estimated by numerical computation.

\end{abstract}


\section{Introduction}
First-passage time (FPT) problems for integrated Markov processes arise both in theoretical and
applied Probability. For instance, in certain stochastic models for the movement of a particle,
its velocity, $Y(t),$ is  modeled as Ornstein-Uhlenbeck (OU) process,
which is indeed suitable to describe the velocity
of a particle immersed in a fluid; as the friction parameter approaches zero, $Y(t)$ becomes
Brownian motion $B_t$ (BM). More generally, the particle velocity $Y(t)$ can be modeled by a diffusion.
Thus, particle position turns out to be the integral of $Y(t),$
and any question about the time at which the particle first reaches a given place leads to the FPT of
integrated $Y(t).$ This kind of investigation  is complicated by the fact that the integral of a Markov process
such as $Y(t),$ is no longer Markovian; however, the two-dimensional process ${\cal Y} (t)= \left (\int _0 ^t Y(s) ds, Y(t) \right )$
is Markovian, so
the FPT of integrated $Y(t)$ can be studied by using Kolmogorov's equations approach. The first apparition in the literature of
${\cal Y} (t),$ with $Y(t)= B_t \ ,$  dates back to the beginning of the twentieth century (see \cite{kolmogorov:anm34}), in modeling a harmonic oscillator
excited by a Gaussian white noise (see  \cite{lachal:crasp97} and references therein).
\par
The study of $\int _0 ^t Y(s) ds$  has interesting  applications
in Biology, in the framework of diffusion models for neural
activity; if one identifies $Y(t)$ with the neuron voltage at time
$t,$ then $\frac 1 t \int _0 ^t Y(s) ds $  represents the time
average of the neural voltage in the interval $[0,t].$ Moreover,
integrated Brownian motion arises naturally in stochastic models
for particle sedimentation in fluids (see \cite{hesse:iamsa05}).
Another application can be found in Queueing Theory, if  $Y(t)$
represents the length of a queue at time $t;$ then $\int _0 ^t
Y(s) ds$ represents the cumulative waiting time experienced by all
the ``users'' till the time $t.$ Furthermore, as an application in
Economy, one can suppose that $Y(t)$ represents the rate of change
of a commodity's price, i.e. the current inflation rate; hence,
the price of the commodity at time $t$ is $X(t)= X(0)+ \int _0 ^t
Y(s) ds .$ Finally, integrated diffusions also play an important
role in connection with the so-called realized stochastic
volatility in Finance (see e.g. \cite{and:01}, \cite{gen:00},
\cite{glot:00}).
\par
FPT problems of integrated BM (namely, when $Y(t)=B_t)$ through
one or two boundaries, attracted the interest of a lot of authors
(see e.g. \cite{ben:13}, \cite{gold:ams71}, \cite{hesse:iamsa05},
\cite{lachal:jap93}, \cite{lachal:aihp91}, \cite{lefebvre:sjam89},
\cite{mckean:jmku63} for single boundary, and
\cite{lachal:crasp97}, \cite{masoliver:phyrev96},
\cite{masoliver:phyrev95} for double boundary); the FPT of
integrated Ornstein-Uhlenbeck process was studied in
\cite{ben:13}, \cite{lefebvre:spa89}. Motivated by these works, our aim is to
extend to integrated Gauss-Markov processes the literature's results
concerning FPT of integrated BM. \par Let $m(t), \ h_1(t), \
h_2(t)$ be $C^1$-functions of $t \ge 0,$ such that $h_2(t) \neq 0$
and $ \rho(t)= h_1(t)/h_2(t) $ is a non-negative and monotonically
increasing function, with $\rho (0)=0.$  \par\noindent If $B(t)=
B_t $ denotes standard Brownian motion (BM), then
\begin{equation} \label{gaussmarkov}
Y(t) = m(t) + h_2(t) B(\rho (t)), \ t \ge 0,
\end{equation}
is a  continuous Gauss-Markov process with mean
$m(t)$ and covariance $c(s,t)= h_1(s) h_2(t), $ for
$0 \le s \le t .$ \par\noindent
Throughout the paper, $Y$ will denote a Gauss-Markov process of the form \eqref{gaussmarkov}, starting from $y =m(0).$
\bigskip

Besides BM,
a noteworthy case of Gauss-Markov process is the Ornstein-Uhlenbeck (OU) process, and in fact
any continuous Gauss-Markov process can be represented in terms of a OU process (see e.g. \cite{ric:smj08}).  \par\noindent
Given a continuous Gauss-Markov process $Y,$ we consider its integrated  process, starting from $X(0):$
\begin{equation} \label{integratedgaussmarkov}
X(t)= X(0) + \int _0^t Y(s) ds .
\end{equation}
For a given boundary $a,$  we study the FPT of $X$
through $a, $ with the conditions that $X(0)=x < a$ and $Y(0)=y,$ that is:
\begin{equation}
\tau _a (x,y)= \inf \{t >0: X(t) =a | X(0) =x, Y(0)=y \};
\end{equation}
moreover, for $b>a$ and $x \in (a,b),$  we also study the first-exit time of $X$ from the interval $(a,b),$
with the conditions that $X(0)=x$ and $Y(0)=y,$ that is:
\begin{equation}
\tau _{a,b} (x,y) = \inf \{t>0: X(t) \notin (a,b) | X(0)=x, Y(0)=y \} .
\end{equation}
In our investigation, an essential role is played by
the representation of $X$ in terms of BM, which was previously obtained by us in \cite{abundo:smj13}. By using this, we avoid to
address the FPT problem by Kolmogorov's equations approach, namely to study the equations associated to the
two-dimensional process $\left (X(t), Y(t) \right );$ many authors
(see the references cited above) followed this analytical approach to study the distribution and the moments of the FPT
of integrated BM, and they obtained explicit solutions, in terms of special functions.
On the contrary, our approach  is based on the properties of Brownian motion and continuous martingales and it has the advantage to be quite simple,
since the problem is reduced to the FPT of a time-changed BM.
Actually, for $Y(0)=y=0$ we present explicit formulae for the density and the moments of the FPT of the integrated Gauss-Markov
process $X,$ both in the one-boundary and
two-boundary case; in particular,
in the two-boundary case, we are able to express the $n$th order moment of the first-exit
time as a series involving only elementary functions. \par

\section{Main Results}
We recall from \cite{abundo:smj13} the following:

\begin{Theorem} \label{proposition1}
Let $Y$ be a Gauss-Markov process of the form \eqref{gaussmarkov};
then $X(t)=x + \int _ 0 ^t Y(s) ds $
is normally distributed with  mean $x+ M(t)$  and variance $\gamma (\rho(t)),$ where
$M(t)= \int _0 ^t m(s) ds,$
$\gamma (t)= \int _0 ^t (R(t) - R(s)) ^2 ds $ and $R (t)= \int _0 ^t h_2(\rho ^{-1} (s))/ \rho '( \rho ^{-1} (s)) ds .$
Moreover, if $\gamma  (+ \infty ) = + \infty, $
then there exists a BM $\widehat B $ such that  $X(t) = x+ M(t) + \widehat B ( \widehat \rho (t) ),$
where $\widehat \rho (t)= \gamma  (\rho(t)).$ Thus, the integrated process $X$ can be represented as
a Gauss-Markov process with respect to a different BM.
\end{Theorem}
\hfill $\Box$

\begin{Remark} {\rm
Notice that, if $\gamma  (+ \infty ) = + \infty, $ though $X$ is represented as a
Gauss-Markov process for a suitable BM $\widehat B, \ X$ is not  Markov
with respect to its natural filtration ${\cal F}_t$ (i.e. the $\sigma-$field generated by $X$ up to time $t)$. In fact,
a Gaussian process $X$ enjoys this property
if and only if
its covariance $K(s,t) = cov (X(s), X(t))$ satisfies
the  condition (see e.g. \cite{doob:49}, \cite{marcus:mar06}, \cite{mehr:65})
$ K(u,t) = \frac {K(u,s)K(s,t)} {K(s,s)}, \ u \le s \le t  .$
Really, if $X$ is e.g. integrated BM with $y=0, \ x=0$
(that is, $X(t)= \int _0 ^t B_s ds),$
one has $K(s,t)  =cov \big ( \int _0 ^s B_u du, \int _0 ^t B_u du \big ) = \frac {s^2} 6 (3t -s)$ (see e.g. \cite{ross:pmodel10}, pg. $654$
or \cite{kle:stoch05}, pg. $105),$ and so the above condition does not hold.
On the other hand,
the two-dimensional process $ \big (\int _0 ^s B_u du, \int _0 ^t B_u du \big )$ has not the same joint distribution as
$\big (\widehat B ( \widehat \rho(s)) , \widehat B ( \widehat \rho(t)) \big ),$ because $cov \big ( \widehat B (\widehat \rho (s)),
\widehat B ( \widehat \rho (t) ) \big )= E [\widehat B (\widehat \rho (s)) \cdot
\widehat B ( \widehat \rho (t) ) ] = \widehat \rho(s) = s^3 /3,$ for $s\le t$ (see Example 1 below), which is different from
$K (s,t) .$  However, the process $(X,B)$ is Markov, and the marginal distributions of the random vector
$ \big (X(s), X(t) \big )$ are equal to the distributions of $\widehat B ( \widehat \rho(s)) $ and $\widehat B ( \widehat\rho(t)), $ respectively; this is enough for the FPT problems we aim to investigate.
}
\end{Remark}

\begin{Remark}  {\rm If $ \gamma  (+ \infty ) = + \infty ,$ and we consider
the time average of $Y$ in the interval $[0,T],$ i.e. $ \overline Y_T =   \frac 1 T  \left ( \int _0 ^T   Y(s)  ds \right ) ,$
by  Theorem \ref{proposition1} we get
$\overline Y_T =  \frac 1 T \left [M(T) + \widehat B (\widehat \rho (T)) \right ],$ namely,
$\overline Y_T $ is normally distributed with mean $(M(T))/T$ and variance $\widehat \rho (T) / T^2 .$  In particular,
if $Y$ is BM, starting from $y$ (that is, $m(t) \equiv y, \ h_2(t) \equiv 1, \ \rho (t) =t),$
 one obtains $\overline Y_T \sim {\cal N } (y,T/3)$ (see  Example 1 below and \cite{abundo:ija08}).
}
\end{Remark}


\bigskip

\noindent {\bf Example 1} (integrated Brownian motion) \par\noindent
Let be $Y(t)=y+ B_t,$  then $ m(t)=y, \ h_1(t) = t,  \ h_2 (t)=1$ and $\rho(t)=t.$ Moreover,
$R(t)=\int _0^t ds =t$ and $\gamma  (t) = \int _0 ^t (t-s)^2 ds = t^3/3 .$
Thus, $\widehat \rho (t) = t^3/3 , \ \gamma  (+ \infty ) = + \infty,$  and so there exists a BM $\widehat B$ such that
$X(t) = x + yt + \widehat B (t^3/ 3)$
(see \cite{abundo:ija08}).
 \par
\hfill $\Box$
\bigskip

\noindent {\bf Example 2} (integrated O-U process)\par\noindent
Let $Y(t)$ be the solution of the SDE (Langevin equation):
$$ dY(t)= - \mu (Y(t)-\beta)dt + \sigma dB_t, \ Y(0)=y,$$
where $\mu , \sigma >0$ and $\beta \in \mathbb{R}.$
The explicit solution is (see e.g.  \cite{abundo:stapro12}):
\begin{equation} \label{ouexplicit}
Y(t)=  \beta + e^{ - \mu t } [ y-\beta +  \widetilde B(\rho (t)],
\end{equation}
where $\widetilde B$ is Brownian motion and $\rho (t)= \frac {\sigma ^2} {2 \mu} \left (e ^{2 \mu
t } -1 \right ) .$ Thus, $Y$ is a Gauss-Markov process with
$m(t)=  \beta + e^{ - \mu t } (y- \beta ), \ h_1(t)= \frac {\sigma ^2} {2 \mu} \left (e ^{\mu
t } -e^{- \mu t } \right ), \
h_2(t)= e^{ - \mu t }$ and $c(s,t)= h_1(s) h_2(t).$ By calculation, we obtain:
\begin{equation} \label{M(t)OU}
M(t)= \int _0 ^t \left ( \beta + e^{ - \mu s } (y- \beta ) \right ) \ ds= \beta t + \frac {(y - \beta)} \mu  \left (1 - e^{ - \mu t } \right ) ,
\end{equation}
\begin{equation}
R(t) = \int _0 ^t  e^{- \mu \rho ^{-1} (s)} ( \rho ^{-1})'(s) ds = \frac {1- e ^{- \mu \rho ^{-1} (t) }  } { \mu} ,
\end{equation}
\begin{equation}
\rho ^{-1} (s)= \frac 1 { 2 \mu} \ln \left (1 + \frac{ 2 \mu} {\sigma ^2} s \right ) ,
\end{equation}
$$ \gamma (t)= \frac 1 { \mu ^2} \int _0 ^t \left ( e^{- \mu \rho ^{-1} (t) } - e^{ - \mu \rho ^{-1} (s)} \right ) ^2 ds =
\frac 1 { \mu ^2} \int _0 ^t \left ( \frac 1 { \sqrt { 1+ 2 \mu t / \sigma ^2 } } - \frac 1 { \sqrt { 1+ 2 \mu s / \sigma ^2 } } \right ) ^2 ds $$
\begin{equation}
=  \frac {\sigma ^2 t } {\mu ^2 (\sigma ^2 + 2 \mu t ) } - \frac {2 \sigma ^2 } {\mu ^3 \sqrt { 1+ 2 \mu t / \sigma ^2 } }
\left ( \sqrt { 1+ 2 \mu t / \sigma ^2 } -1 \right ) + \frac { \sigma ^2 } {2 \mu ^3 } \ln \left ( 1+ 2 \mu t / \sigma ^2 \right ) .
\end{equation}
Then, by Theorem \ref{proposition1}, we get that
$X(t)= x+ \int _0 ^t Y(s) ds $ is normally distributed with mean
$x+ M(t)$
and variance $\widehat \rho (t) = \gamma  (\rho (t))  .$ Moreover, as easily seen,
$\lim _ { t \rightarrow + \infty} \gamma (t) = + \infty ,$ so
there exists a BM $\widehat B $ such that $X(t)= x+ M(t) + \widehat B \left ( \widehat \rho (t) \right ).$
\hfill $\Box$
\bigskip

\noindent
Notice that in both Example 1 and 2 it holds $\rho ( + \infty )= + \infty,$ so the condition $\gamma  ( + \infty) = + \infty$
is equivalent to $\widehat \rho ( + \infty ) = + \infty .$
\bigskip

\noindent{\bf Example 3} (integrated Brownian bridge) \par\noindent
For $T>0$ and  $\alpha, \beta \in \mathbb{R},$ let $Y(t)$ be the solution of the SDE:
$$ dY(t)= \frac {\beta - Y(t)} {T-t } \ dt +  dB_t, \ 0 \le t \le T, \ Y(0)= y = \alpha .$$
This is a transformed BM with fixed values at each end of the interval $[0,T], \ Y(0)=y= \alpha $ and $Y(T)=\beta.$
The explicit solution is (see e.g. \cite{revuzyor:con91}):
$$ Y(t) =\alpha \left ( 1 - t/T \right ) + \beta t /T + (T-t) \int _0 ^t \frac 1 { T-s} dB(s) $$
\begin{equation} \label{explicitBbridge}
= \alpha \left ( 1 - t/T \right ) + \beta t /T + (T-t) \widetilde B \left ( \frac { t} {T(T-t) } \right )  , \ 0 \le t \le T ,
\end{equation}
where $\widetilde B$ is BM.
So, for $0 \le t \le T, \ Y$ is a Gauss-Markov process with:
$$m(t)= \alpha \left ( 1 - t/T \right ) + \beta t /T, \  h_1(t) = t/T , \ h_2(t) =T-t, \ \rho (t) = \frac t { T(T-t)} , \ c(s,t)= h_1(s)h_2(t) .$$
Notice that now $\rho (t)$ is defined only in $[0,T).$
By calculation, we obtain:
\begin{equation} \label{M(t)brownianbridge}
M(t)= \alpha t + \frac {\beta - \alpha } {2T }  t^2, \
R(t) = \frac {T^3t (2+Tt) } {2(1+Tt)^2}  ,
\end{equation}
$$ \rho ^{-1} (s)= \frac {T^2 s } {1+Ts } , \
 \gamma (t)= \int _0^t \left ( \frac {T^3t (2+Tt) } {2(1+Tt)^2}- \frac {T^3s (2+Ts) } {2(1+Ts)^2} \right )^2   ds.$$
Then, by Theorem \ref{proposition1}, we get that $X(t)=
x+ \int _0 ^t Y(s) ds $ is normally distributed with mean $x+ M(t)$ and
variance $\widehat \rho (t)= \gamma  (\rho (t))  .$ As easily seen,
$\lim _ {t \rightarrow T^-} \rho (t) = + \infty ;$ moreover,
by a straightforward, but
boring calculation, we get that  $\lim _ { t
\rightarrow T^-} \widehat \rho (t)=  \gamma_1 ( + \infty )= + \infty ,$ so there exists a
BM $\widehat B $ such that $X(t)= x+ M(t) + \widehat B \left (
\widehat \rho (t) \right ) , \ t \in [0,T].$ \par\noindent
\hfill $\Box$
\bigskip

\noindent {\bf Example 4} (the integral of a generalized Gauss-Markov process) \par\noindent
Let us consider the diffusion $Y(t)$ which is the solution of the SDE:
$$ dY(t)= m'(t) dt + \sigma (Y(t)) dB_t , \ Y(0) = y, $$
where  $\sigma (y)>0$ is a smooth deterministic function. In this Example, we denote by $\rho (t)$
the quadratic variation of $Y(t),$ that is,
$\rho (t):= \langle Y \rangle _t = \int _0 ^t \sigma ^2 (Y(s)) ds , $ and suppose that
 $\rho ( + \infty)= + \infty .$
By using the Dambis, Dubins-Schwarz Theorem (see e.g. \cite{revuzyor:con91}), it follows that
$Y(t)= m(t)+ \widehat B( \rho (t)), \ t \ge 0  \ (m(0)=y) ,$ where $\widehat B$ is BM;
here, $\rho(t)$ is increasing, but not necessarily deterministic, namely
it can be a random function.
For this reason, we call $Y$ a generalized Gauss-Markov process.
Denote by $A$ the ``inverse'' of the random function $\rho ,$ that is,
$ A(t) = \inf \{ s >0 : \rho (s) >t \};$ since $\rho(t)$ admits derivative and $\rho'(t)=\sigma ^2 (Y(t)) > 0,$ also $A'(t)$
exists and
$ A'(t)= \frac 1 { \sigma ^2 (Y(A(t))};$
we focus on the case when
there exist  deterministic continuous  functions $\alpha (t), \ \beta (t) $ (with  $\alpha(0)= \beta (0)=0)$ and
$\alpha _1(t), \ \beta _1 (t),$ such that, for every $t \ge 0:$
$$ \alpha (t), \  \beta (t) \ {\rm are \ increasing ,} \ \alpha (t) \le \rho (t) \le \beta (t) , \ {\rm and} \
\alpha _1 (t) < A'(t) < \beta _1  (t).$$
Since $\rho (t)$ is not, in general, deterministic, we cannot obtain exactly the distribution of  $ \int _0 ^t Y(s) ds,$ however we are able
to find lower and upper bounds to it. In fact, we have:
$$ \int _0^t Y(s) ds = \int _0 ^t m(s) ds + \int _0 ^t \widehat B(\rho (s)) ds=
\int _0 ^t m(s) ds + \int _0 ^ {\rho (t)} \widehat B(v) A' (v) dv .$$
By using the arguments leading to the proof of Theorem \ref{proposition1},
(see \cite{abundo:smj13} for more details), we conclude that, for fixed $t$ the law of $\int _0 ^t Y(s)ds ,$ conditional to $ \rho (t),$
is normal with mean $ M(t)= \int _0 ^t m(s) ds$ and variance $\gamma  (\rho(t)) ,$
which is bounded between $\gamma (\alpha (t))$ and $\gamma (\beta (t)).$
Here, $\gamma  (t)= \int _0 ^t (R(t)- R(s))^2 ds,$ where $R(t) = \int _0 ^t A'(s) ds$
is bounded between $\int _0 ^t \alpha _1 (s) ds$ and $\int _0 ^t \beta _1 (s) ds .$
The closer $\alpha (t)$ to $\beta (t)$ and $\alpha _1 (t)$ to $\beta _1 (t),$ the better the approximation above;
for instance, if $\sigma (y)= 1+ \epsilon \cos ^2 (y), \  \epsilon >0 ,$ we have
$\rho(t) = \int _0 ^t (1+  \epsilon \cos ^2 (Y(s)))^2 ds$ and so
$\alpha (t)= t, \ \beta (t)= (1+ \epsilon)^2 t , \ \alpha _1(t) = 1 / (1+ \epsilon )^2, \ \beta _1(t)= 1. $
The smaller  $\epsilon,$ the closer $\gamma (\alpha (t))$ to $\gamma (\beta (t)).$ \par
\hfill $\Box$

\bigskip

\bigskip

In the sequel, we suppose that all the assumptions of Theorem \ref{proposition1} hold,
and $ \gamma  (+ \infty )= + \infty;$
we limit ourselves
to consider the special case when $m(t)$ is a constant (that is, $m(t) \equiv  Y(0)=y, \ \forall t ) ,$
thus $Y(t) = y + h_2(t) B( \rho (t))$ and
$X(t)= x + yt + \int _0 ^t h_2(s) B( \rho (s)) ds .$
Our aim is to investigate the FPT problem of $X,$ for one or two boundaries. One approach to the FPT problem of $X$ consists
in considering
the two-dimensional process $(X(t), Y(t))$ given by:
$$ \begin{cases}
X(t)= x + \int _0 ^t Y(s)ds \\
Y(t)=y + h_2(t) B( \rho (t)) dt \ ,
\end{cases}
$$
or, in differential form:
$$ \begin{cases}
dX(t)= Y(t)dt \\
dY(t)=h_2'(t) B( \rho (t)) dt + h_2(t) \sqrt { \rho ' (t)} d B_t \ ,
\end{cases}
$$
and to study the associated Kolmogorov's equations. \par\noindent
Many authors
(see e.g. \cite{gork:ams75}, \cite{lachal:crasp97},
 \cite{lachal:jap93}, \cite{lachal:aihp91}, \cite{lefebvre:sjam89})
followed this way in the case of
integrated BM, namely for $Y(t)=y+B_t \ .$  In fact, for $ \tau = \tau _a$ or $\tau = \tau _{a,b},$ the law of the couple
$(\tau  (x,y), B_ { \tau  (x,y)})$ was investigated. Let us denote by ${\cal G}$  the generator of
$(X, B) ,$  that is:
$$ {\cal G} f(x,y)=  \frac {\partial f } { \partial x} \cdot y + \frac 1 2 \frac {\partial ^2 f } { \partial y ^2}  \  , \ f \in C^2 ;$$
if one considers, for instance, the one boundary case, then the Laplace transform of \par\noindent
$ \left (\tau _a (x,y), B_ { \tau _a (x,y)} \right ),$ defined
for $x \le a , \ y \in \mathbb{R} ,$
by $u(\lambda, \nu):= E  \left [ \exp \left (- \lambda \tau _a (x,y) - \nu B_ {\tau _a (x,y) } \right ) \right ]$
 \par\noindent
 $(\lambda, \ \nu \ge 0) ,$ is the solution of the problem with boundary conditions:
\begin{equation}
\begin{cases}
{\cal G} u(x,y) = \lambda u(x,y), \ x \le a , \ y \in \mathbb{R} \\
u(a^-, y)= e^{ - \nu y}, \ y \ge 0 \\
u(a^+, y)= e^{  \nu y}, \ y < 0
\end{cases}
\end{equation}
(see e.g. \cite{lachal:aihp91}, Lemma 3, or ref. [4], [5], [7], therein).
Moreover,  for $n=1, 2, \dots $ the $n$th order moments $T_n(x,y)= E( \tau _a ^n (x,y)) $
are solutions to the equations ${\cal G}T_n = - n T_ {n-1} \ (T_0  \equiv 1),$ subjected to certain boundary conditions; however, these
boundary value problems are not well-posed (see \cite{hesse:iamsa05},
where some numerical methods to estimate  $T_n$ were also considered).
\par
Notice that, in the case of integrated BM, explicit, rather complicated  formulae for the joint distribution of
$ \left (\tau_a (x,y) , B_ {\tau_a (x,y)} \right )$ (and therefore for the density of $\tau_a (x,y))$ were found in \cite{gold:ams71},
\cite{lachal:jap93}, \cite{mckean:jmku63}).
In order to avoid not convenient formulae,
we propose an alternative approach, based on the representation of the integrated process $X$ as a Gauss-Markov process,
with respect to  the BM $\widehat B  $ (see Theorem \ref{proposition1}); this way works very simply, almost in the case when
$y=0.$
Thus, in the following, we suppose that $Y(t) = y + h_2(t) B( \rho (t))$ and
$\gamma  (+ \infty ) = + \infty, $ so the integrated  process is
of the form $X(t)= x + yt + \widehat B (\widehat \rho (t)),$ where $\widehat \rho (t)= \gamma  ( \rho (t))$ and
$\widehat B$ is a suitable BM.
Notice however, that the integrated OU process and  the integrated Brownian bridge belong to this class only if
 $y= \beta $ (see \eqref{M(t)OU}), and $\alpha = \beta =y$ (see \eqref{M(t)brownianbridge}), respectively.

\subsection{FPT through one boundary}
Under the previous assumptions, let  $a$ be a fixed constant boundary;
for $x<a$ and $y \in \mathbb{R},$  the FPT of $X$ through $a$ can be written as follows:
\begin{equation}
\tau _a (x, y) = \inf \{t >0: x+ yt+ \widehat B (\widehat \rho (t))  =a  \} .
\end{equation}
Thus, if we set $\widehat \tau _a (x, y) = \widehat \rho ( \tau _a (x, y)),$ we get:
\begin{equation} \label{hattau}
\widehat \tau _a (x, y) = \inf \{t >0:  \widehat B_t  = h(t)  \},
\end{equation}
where $h(t)= a-x -y   \widehat \rho \ ^{-1} (t) ,$
and so we reduce to consider the FPT of BM through a curved boundary. Since, for $x <a$ and $y \ge 0$ the
function $h(t)$ is not increasing, we are able to conclude that $\tau _a (x,y)$ is finite with probability one, if $y \ge 0.$
In fact, as it is well-known, the FPT of BM $\widehat B_t$ through the constant barrier $h(0)=a-x,$ say
$\bar \tau (x),$ is finite with probability one;
then, if $y \ge 0,$  from $h(t) \le h(0)$ we get that $\widehat \tau _a (x, y) \le \bar \tau (x)$ and therefore
also $\widehat \tau _a (x, y)$ is finite with probability one. Finally, if $y \ge 0,$ we obtain that
$P( \tau _a (x,y) < + \infty ) =1,$
because $\tau _a (x,y) = \widehat \rho ^ {-1} ( \widehat \tau _a (x,y)) \le \widehat \rho ^ {-1} ( \bar \tau _a (x)).$
Note, however, that this argument does not work for $y <0.$ \par\noindent
A more difficult problem is to find the distribution of $\widehat \tau _a (x, y),$ and then that of $ \tau _a (x, y) .$
However, if $h(t)$ is either convex or concave, then
lower and upper bounds to the distribution of
$\widehat \tau _a (x, y)$ can be obtained by considering a ``polygonal approximation'' of $h(t)$ by means of a piecewise-linear function
(see e.g. \cite{abundo:ijam12}, \cite{abundo:ric01}), but
in general, it is not possible to find the distribution of $\widehat \tau _a (x, y)$ exactly.
\bigskip

\begin{Remark} {\rm
Actually, it is possible to find explicitly the density of the FPT of $X$ through certain moving boundaries.
Indeed, denote by ${\cal V}$  the family of continuous functions  which consists of curved boundaries $v=v(t), \ t \ge 0 , \ v(0)>0,$
for which the FPT-density of BM through $v$ is explicitly known; this family
includes linear boundaries $v(t)=at + b$ (see \cite{bache:41}), quadratic boundaries $v(t)=a - b t^2$
(see e.g. \cite{ali:10}, \cite{groe:ptr89}, \cite{salm:aap88}), square root boundaries
$v(t)= a\sqrt { t+b},$ and $ v(t)= a \sqrt { (1+bt) (1+ct)}$ (see e.g. \cite{ali:10}, \cite{bre:67}, \cite{sato:77}),
and the
so-called Daniels boundary $v(t)= \delta - \frac t {2 \delta} \log \left ( \frac {k_1} 2 + \sqrt { \frac { k_1 ^2} 4 + k_2
e^{ - 4 \delta ^2 /t } } \right ) $
(see \cite{dan:ast82}, \cite{dan:jap69}). For a boundary $v \in {\cal V},$
denote by $\widehat f_v (t |x)$ the FPT-density of BM starting from $x < v(0)$ through the boundary $v;$
if $S(t)= v( \widehat \rho (t)) + yt ,$
then one can easily find the density of the FPT of $X$ through $S,$ with the condition that $x < S(0)=v(0).$ In fact, if
 $\tau _S (x,y) = \inf \{ t>0: X(t) = S(t) |X(0)=x, Y(0)=y \},$ one gets
$\tau _S (x,y) = \inf \{ t>0: x +ty + \widehat B ( \widehat \rho (t)) = S(t) \};$ then,
 $\widehat \tau _v (x,y) := \widehat \rho (\tau _S (x,y)) = \inf \{t>0: x + \widehat B (t) = v(t) \}$ has density $\widehat f_v$ and so
the density of $ \tau _S (x,y)$ turns out to be
\begin{equation} \label{fptdenoverv}
f_S (t|x)= \widehat f_v (\widehat \rho (t)|x) \widehat \rho '(t) .
\end{equation}

For instance, if $X$ is integrated BM $(\widehat \rho (t)= t^3/3)$, and we consider
the cubic boundary $S(t)=a+ty +b t^3$ \
$(a >0, b <0 ),$  it results $S(t)= v( \widehat \rho (t)) + yt,$ with $v(t)= a+ 3bt $
and so, for $x < a, \ \widehat \tau _v(x,y)$
is the FPT of BM starting from $x$ through the linear boundary $a + 3bt.$ Thus, $\widehat \tau _v (x,y)$ has the inverse Gaussian
density $\widehat f _v(t|x)= \frac {a-x } {\sqrt {2 \pi } \ t^ {3/2} } e^ { - (3bt+a-x)^2 /2t }$ (see e.g. \cite{abundo:ric01});
then, the density of $\tau _S (x,y)$ is obtained by  \eqref{fptdenoverv}.
}
\end{Remark}
\bigskip

Formula \eqref{hattau}, with $y=0,$ allows to find the density of $ \tau _a (x, 0)$ in closed form; in fact,
$\widehat \tau _a (x, 0)$ is the FPT of BM $\widehat B$ through the level $a-x>0,$ and so its density is:
\begin{equation} \label{hatfa}
\widehat f _a(t|x) := \frac d {dt} P( \widehat \tau _a (x, 0) \le t )= \frac {a-x } {\sqrt {2 \pi } \ t^ {3/2} } e^ { - (a-x)^2 /2t },
\end{equation}
from which the density of $\tau _a (x, 0)= \widehat \rho ^ {-1} ( \widehat \tau _a (x,0)$ follows:
\begin{equation} \label{densityoftaua}
 f _a(t|x) := \frac d {dt} P( \tau _a (x, 0) \le t )= \widehat f _a (\widehat \rho (t) |x ) \widehat \rho ' (t)
 = \frac {(a-x) \ \widehat \rho ' (t) } {\sqrt {2 \pi } \  \widehat \rho (t)^ {3/2} } e^ { - (a-x)^2 /2 \widehat \rho (t) }.
\end{equation}
If $X$ is integrated BM, we have $X(t)= x +  \widehat B (\widehat \rho (t)),$ with
$\widehat \rho (t) = t^3/3,$ so we get
(cf. \cite{gold:ams71}):
\begin{equation} \label{IBMdensitythrougha}
 f_a (t|x) = \frac {3^ {3/2}(a-x) } {\sqrt {2 \pi } \ t^ {5/2} } e^ { - 3(a-x)^2 /2t^3  }.
\end{equation}
If $X$ is integrated OU process, the density of $ \tau _a (x,0)$ can be obtained by inserting in \eqref{densityoftaua} the
function $\widehat \rho (t)$ deducible from Example 2, but it takes a more complex form.

\begin{Remark}
{\rm
Formula \eqref{densityoftaua} implies that the $n$th order moment of the FPT, $E( \tau _a ^n (x, 0)),$ is finite
if and only if the function $t^n \widehat \rho '(t)/ \widehat \rho (t) ^{3/2}$ is integrable in $(0, + \infty).$ \par\noindent
Now, let us suppose  that there exists $\alpha >0$ such that $\widehat \rho (t) \sim const \cdot t^ \alpha,$ as $ t \rightarrow + \infty ;$
then, in order that
$E( \tau _a ^n (x, 0)) < \infty , $  it must be $\alpha = 2(n + \delta ),$
for some $\delta  >0.$
For integrated BM,
we have $\alpha =3,$ then for $n=1$  the last condition holds with $\delta = 1/2,$ so we obtain  the
finiteness of $E( \tau _a  (x, 0))$ (notice that the mean FPT of BM through a constant barrier is instead infinite).
Of course, this is not always the case; in fact, if $X$ is  integrated OU process,
we have $\rho (t) \sim const \cdot  e^{2 \mu t}, \ \gamma  (t) \sim const \cdot \ln ( 2 \mu t / \sigma ^2 ), $ as $ t \rightarrow + \infty ,$
and so $ \widehat \rho (t) = \gamma  ( \rho (t)) \sim const \cdot t ,$  as $t \rightarrow + \infty , $ namely $\alpha =1$ and the condition above is not satisfied
with $n=1;$ therefore $E( \tau _a (x, 0)) = + \infty.$ Not even $ E( ( \tau _a (x, 0) )^ {1/2})$ is finite, but $ E((\tau _a (x,0))^{1/4})$ is so.
Notice that the moments of any order of the FPT of (non integrated) OU through a constant barrier are instead finite.
\par
As for the second order moment of the FPT of
integrated BM, instead, we obtain
$E\left [ \left ( \tau _a (x, 0) \right )^2 \right ] = + \infty,$ since the equality $\alpha = 2(n
+ \delta )$ with $\alpha =3$ and $n=2$ is not satisfied, for any
$\delta >0.$ }
\end{Remark}
\par
From \eqref{hatfa} we get that the $n$th order moment of $\tau _a (x,0),$  if it exists finite, is explicitly given by:
$$
E \left [ ( \tau _a (x,0))^n \right ] = E \left [ (\widehat \rho \ ^ {-1} ( \widehat \tau _a (x, 0)))^n \right ]
$$
\begin{equation}
= \int _0 ^ { + \infty }  ( \widehat \rho \ ^ {-1} (t))^n \frac {a-x } {\sqrt {2 \pi } t^ {3/2} } e^ { - (a-x)^2 /2t }  dt .
\end{equation}
For instance, if $X$ is integrated BM,  one has:
$$
E( \tau_a (x, 0)) = E( (3 \ \widehat \tau _a (x, 0))^ {1/3} ) = \int _0 ^ {+ \infty } (3t)^ {1/3} \frac {a-x} {\sqrt {2 \pi } t^ {3/2} }  e^ { - (a-x)^2 /2t }  dt
$$
$$
= \frac { 3 ^ {1/3} (a-x) } {\sqrt {2 \pi } } \int _ 0 ^ { + \infty } \frac 1 { t ^ {7/6} } e^ { - (a-x)^2 /2t }  dt .$$
By the variable's change  $z = 1/t,$ the  integral can be written as:
$$  \int _ 0 ^ { + \infty } \frac 1 { z ^ {5/6} } e^ { - (a-x)^2 z /2 }  dz =
\frac {\Gamma \left ( \frac 1 6 \right ) 2^ {1/6} } {(a-x)^ {1/3} }  \int _0 ^ {+ \infty } \left ( \frac {(a-x)^2 } 2 \right ) ^ {1/6}
\frac 1 { \Gamma \left ( \frac 1 6 \right ) } z ^ {1/6 -1} e^ {- \frac { (a-x)^2 } 2 z } dz $$
$$ = \frac {\Gamma \left ( \frac 1 6 \right ) 2^ {1/6} } {(a-x)^ {1/3} } \ ,
$$
where we have used that the last integral equals one, because the integrand is a Gamma density.
Thus, for integrated BM, we finally obtain:
\begin{equation} \label{meantauay0}
E( \tau _a (x, 0))=
 \left ( \frac 3 2 \right ) ^{1/3} \Gamma \left (\frac 1 6 \right ) \frac {  (a-x)^ {2/3} } {\sqrt { \pi } }  \  .
\end{equation}

Until now we have supposed that the starting point $x < a$ is given and fixed. We can introduce a randomness in the starting point, replacing $X(0)=x$ with a random variable $\eta ,$ having  density $g(x)$ whose support is the interval $(- \infty, a );$ the corresponding FPT problem is particularly relevant in contexts such as neuronal modeling, where the reset value of the membrane potential is usually unknown (see e.g. \cite{lansky:89}). In fact,
the quantity of interest becomes now the unconditional FPT through the boundary $a,$ that is,
$ \inf \{ t >0: X (t) = a | Y(0)=y \} ;$ in particular, if $X$ is integrated BM and $y=0,$ one gets from \eqref{meantauay0} that
the average FPT through the boundary $a, $ over all initial positions $\eta <a,$ is:
\begin{equation} \label{meantauy0random}
\overline T _a = \int _ { - \infty} ^a E( \tau _a (x,0)) g(x) dx = \left ( \frac 3 2 \right ) ^{1/3} \frac {\Gamma \left (\frac 1 6 \right )} {\sqrt { \pi } } \int _ {- \infty } ^a  (a-x)^ {2/3} g(x) dx .
\end{equation}
For instance, suppose that  $a- \eta$ has Gamma distribution with parameters $ \alpha , \  \lambda >0 ,$ namely, $\eta $ has density
$$ g(x)= \frac {\lambda ^ \alpha } { \Gamma ( \alpha )} e^ {- \lambda (a-x)} (a-x) ^ {\alpha -1} \cdot
\mathbb{I} _ { (- \infty , a ) } (x) .$$
Then, by the change of variable $z=a-x $ one obtains that the above integral is nothing but $E \left ( Z ^ {2/3} \right ),$ where
$Z$ is a random variable with the same distribution of $a- \eta ;$ then, recalling the expressions of the moments of the Gamma distribution, one obtains
$E \left ( Z ^ {2/3} \right )= \frac {\Gamma ( \alpha + \frac 2 3 ) } { \lambda ^ {2/3} \Gamma ( \alpha )} . $
Finally, by inserting this quantity in \eqref{meantauy0random}, it follows that:
$$ \overline T _a = \frac { \left ( \frac 3 {2 \lambda ^2 }\right ) ^ {1/3}   }  {\sqrt \pi  } \cdot \frac { \Gamma \left ( \frac 1 6 \right )
\Gamma \left ( \alpha + \frac 2 3 \right ) }  { \Gamma ( \alpha ) }.$$

\begin{Remark}
{\rm
For $y= Y(0)=0,$ we have considered the FPT of $X$ through the boundary $a$ from ``below'', with the condition $x= X(0)<a;$
if one considers the FPT of $X$ through the barrier $a$ from ``above'', with the condition $ X(0)>a$
(namely, $\inf \{ t>0: X(t) \le a | X(0) =x, Y(0)=0 \}),$
then in all formulae $a-x$  has to be replaced with $x-a.$ More generally, if one considers the first hitting time of $X$ to $a$ (from above or below), $a-x$ must be replaced by $|a-x|.$
\bigskip
}
\end{Remark}

\subsection{FPT in the two-boundary case: first exit time from an interval}
Assume, as always, that $\gamma (+ \infty ) = + \infty ;$ for $x \in (a,b)$ and $y \in \mathbb{R},$
the first-exit time of $X$
from the interval $(a,b)$ is:
\begin{equation}
\tau _{a,b} (x,y) = \inf \{t >0: x+ yt+ \widehat B (\widehat \rho (t)) \notin (a,b) \} .
\end{equation}
Set $\widehat \tau _{a,b} (x, y) = \widehat \rho ( \tau _{a,b} (x, y)),$ then:
\begin{equation}
\widehat \tau _{a,b} (x, y) = \inf \{t >0: x+ \widehat B_t  \le a -y   \widehat \rho  ^ {-1} (t) \ \  {\rm or } \  \
 x + \widehat B_t  \ge b -y   \widehat \rho  ^ {-1} (t) \} .
\end{equation}
If $\widehat \tau _{a,b} (x, y)$ is finite with probability one, also $ \tau _{a,b} (x, y)$ is so.
In the sequel, we will focus on the case when  $y=0,$ namely we will consider $\tau _{a,b} (x, 0) = \widehat \rho ^  {-1} (\widehat \tau _{a,b} (x, 0) ),$
where
$ \widehat \tau _{a,b} (x, 0) = \inf \{t >0: x + \widehat B_t  \notin (a, b) \} ;$
as it is well-known, $ \widehat \tau _{a,b} (x, 0)$ is finite with probability one and its moments are
solutions of Darling and Siegert's equations (see \cite{darling:ams53}). \par
First, we will find sufficient conditions so that the moments of $ \tau _{a,b} (x, 0)$ are finite; then,
we will carry on explicit computations of them, in the case of integrated BM. \par\noindent
\begin{Proposition}
If $\widehat \rho $ is convex, then $E \left ( \tau _{a,b} (x,0) \right ) < \infty ;$ moreover, if there exist constants $c, \ \delta >0,$
such that $0 \le \widehat \rho  ^ {-1} (t) \le c \cdot  t ^ \delta, $ then
$ E \left ( \tau _{a,b} (x, 0) \right ) ^n < \infty ,$ for any integer $n.$
\end{Proposition}
\bigskip

\noindent {\it Proof.} \
If $\widehat \rho $ is convex, then $\widehat \rho ^ {-1} $ is concave,
and the finiteness of $E \left ( \tau _{a,b} (x,0) \right )$ follows by Jensen's inequality written for concave functions. Next,
denote by $\widehat f _ {- \alpha, \alpha } (t |x) $ the density of the first-exit time
of $x+ \widehat B_t$ from the interval $(- \alpha, \alpha ), \ \alpha >0;$ we
recall from \cite{darling:ams53} that
the Laplace transform of $\widehat f _ {- \alpha, \alpha } (t |x),$ namely,
$ \int _ 0 ^ { + \infty } e^{- \theta t} \widehat f _ {- \alpha, \alpha } (t |x) dt $ is:

\begin{equation} \label{laplacefpttwobarrier}
{\cal L} \left [ \widehat f _ { - \alpha, \alpha } \right ] ( \theta |x )=
 \frac {\cosh (\sqrt
{ 2 \theta } x  )} { \cosh (\sqrt { 2 \theta } \alpha ) } \ , \
-\alpha < x < \alpha, \ \theta \ge 0  .
\end{equation}
By inverting this Laplace transform, one obtains (see \cite{darling:ams53}):
\begin{equation} \label{fptdensitytwobarrier}
\widehat f_{- \alpha, \alpha } (t| x  )= \frac \pi {\alpha ^2 } \sum _ {
k=0} ^ \infty (-1)^k \left (k + \frac 1 2 \right ) \cos \left [
\left (k + \frac 1 2 \right ) \frac { \pi x } \alpha \right ] \exp
\left [ -  \left (k + \frac 1 2 \right ) ^2 \frac {x^2 t } {2
\alpha ^2 } \right ].
\end{equation}
The case of an interval $(a, b), \ b >a ,$ is reduced
to the previous one; in fact, as easily
seen, if $ \alpha = (b -a )/2$ one has:
$$
\widehat f_{a, b } (t| x  )= \widehat f_{- \alpha, \alpha } \left (t| x - \frac
{a+b} 2 \right ).
$$
Of course, the density of $\tau _{a, b} (x,0)$ turns out to be $\widehat f _{a,b} ( \widehat \rho (t)|x ) \widehat \rho '(t).$
For the sake of simplicity, we take $a=-\alpha, \ b = \alpha, \ \alpha >0;$ then, for $x \in (- \alpha, \alpha)$ and an integer $n:$
\begin{equation} \label{tauIBM}
E \left [ \left ( \tau _{a,b} (x, 0) \right )^n \right ]= E \left [ \left ( \tau _{- \alpha, \alpha} (x, 0) \right ) ^n \right ] =
E \left [ \left ( \widehat \rho ^ {-1} ( \widehat \tau _{- \alpha, \alpha} (x, 0) \right ) ^n  \right ]
= \sum _ {k=0} ^ \infty A_k (x),
\end{equation}
where
\begin{equation} \label{Ak}
A_k(x)= \frac \pi { \alpha ^2 } (-1)^k \left (k + \frac 1 2 \right ) \cos \left ( \left (k + \frac 1 2 \right ) \frac { \pi x } \alpha \right )
\int _0 ^ { + \infty } e^ { - (k+ 1/2)^2 \pi ^2 t / 2 \alpha ^2 } \left ( \widehat \rho ^{-1} (t) \right ) ^n dt .
\end{equation}
The integral can be written as:
$$ \frac {2 \alpha ^2 } {\pi ^2 (k +1/2) ^2 } E \left ( \widehat \rho ^ {-1} (Z_k) \right ) ^n,$$
where $Z_k$ is a random variable with exponential density of parameter $\lambda _k= (k+1/2)^2 \pi ^2 / 2 \alpha ^2 ;$
so:
$$ A_k (x)= (-1)^k  \cos \left ( \left (k + \frac 1 2 \right ) \frac { \pi x } \alpha \right )
\frac {2  } {\pi  (k +1/2)  } E \left ( \widehat \rho ^ {-1} (Z_k) \right ) ^n.$$
Recalling that $E[ (Z_k)^{n \delta } ]= \frac { \Gamma (1+n \delta) } {(\lambda _k) ^{n \delta }  },$
by the hypotheses we get $E \left ( (\widehat \rho ^ {-1} (Z_k))^n \right )  \le c ^n E [ (Z_k) ^ {n \delta }] =
const  \cdot \frac {\Gamma (1 + n \delta  ) } {(k+ 1/2)^{2  n \delta }  } ;$
thus:
$$ | A_k(x)| \le  \frac {const '} { (k+1/2)^ {1+2n \delta } } ,$$
from which it follows that the series $\sum _k A_k(x)$ is absolutely convergent for every $x \in ( - \alpha, \alpha),$ and therefore
$E \left [ \left ( \tau _{- \alpha, \alpha } (x, 0) \right )^n \right ] < + \infty .$
The finiteness of $E \left [ \left (  \tau _{a, b} (x, 0) \right )^n  \right ]$ in the general case is easily obtained.
\par \hfill  $\Box$
\bigskip

\begin{Remark}
{ \rm The condition
$0 \le \widehat \rho ^ {-1} (t) \le c \cdot t ^ \delta $ is satisfied e.g. for integrated BM, since  $ \widehat \rho ^ {-1} (t) = 3^ {1/3} t^ {1/3}$
(see Example 1), and for integrated OU process, because from the expression of $\widehat \rho (t)$ deducible from Example 2, it can be shown that
$c_1 t \le  \widehat \rho (t) \le c_2 t $ for suitable $c_1, c_2 >0$ which depend on $\mu$ and $\sigma ,$  and therefore
$\frac 1 {c_2} t \le \widehat \rho ^ {-1} (t) \le \frac 1 {c_1} t .$ }
\end{Remark}

Now, we carry on explicit computations of $E \left [ \tau _{a, b} (x, 0) \right ] $ and
$E \left [ \left (  \tau _{a, b} (x, 0) \right )^2 \right ] ,$
in the case of integrated BM.
Inserting  $ \widehat \rho (t) = t^3/3, \ (\widehat \rho ^ {-1} (y)= (3y) ^ {1/3}), $  and $n=1, 2$ in \eqref{tauIBM}, \eqref{Ak},
after some calculations we obtain:
\begin{equation} \label{explictauIBM}
E\left [ \tau _{a, b} (x, 0) \right ] = \frac {3^ {1/3} 2^{7/3} \Gamma (\frac 4 3 ) (b-a) ^ {2/3}  } {\pi ^ {5/3} }
\sum _ {k=0} ^ \infty (-1)^k \frac 1 {(2k+1) ^ {5/3} }
\cos \left [ \frac { \pi (2k+1)} {b-a}  \left (x - \frac {a+b} 2 \right )  \right ] .
\end{equation}

\begin{equation} \label{explictauquadIBM}
E \left [ \left ( \tau _{a, b} (x, 0) \right ) ^2 \right ] = \frac {12 (b-a) ^ 4  } {\pi ^ 4 }
\sum _ {k=0} ^ \infty (-1)^k \frac 1 {(2k+1) ^ 4 }
\cos \left [ \frac { \pi (2k+1)} {b-a}  \left (x - \frac {a+b} 2 \right )  \right ] .
\end{equation}
Notice that it is arduous enough to express the sums of the Fourier-like series above in terms of elementary functions of $x \in (a,b),$ and then
to obtain the moments of $\tau _{a, b} (x, 0) $ in a simple closed form; actually, by
using the Kolmogorov's equations approach,
in \cite{masoliver:phyrev96}, \cite{masoliver:phyrev95},  it was obtained
a formula for $E( \tau _ {a,b} (x,0))$  in terms of hypergeometric functions.
This kind of difficulty does not arise, for instance, in the case of (non-integrated)
BM; in fact, by using formula  \eqref{tauIBM}
with $\widehat \rho (t)=t$ and $n=1,$
 one obtains:
$$ E\left [ \tau _{- \alpha, \alpha} (x) \right ]= \frac {32 \alpha ^2 } { \pi ^3 } \sum _ {k=0} ^ \infty (-1)^k \frac 1 {(2k+1)^3 }
\cos \left [ (2k+1) \frac \pi {2 \alpha} x \right ] ;$$
on the other hand, the well-known formula for the
mean first-exit time of BM from the
interval $(- \alpha, \alpha),$ provides
that the sum of the series must be $\alpha ^2 - x^2.$
\par
However, \eqref{explictauIBM} and \eqref{explictauquadIBM} turn out to be very convenient to estimate the first two moments of $\tau _{a, b} (x, 0)$ for integrated BM; in fact
the two series  converge fast enough, so to obtain ``good'' estimates of the moments,  it suffices to consider a few terms of them. As for
$E\left [ \tau _{a, b} (x, 0) \right ],$ it
appears to be fitted
very well by the square root of a quadratic function.
In the Figure 1, for integrated BM, we compare the graphs of $E( \tau _{a, b} (x, 0)),$ calculated by replacing the series in \eqref{explictauIBM} with a finite summation over the first $20$ addends,
and that of  $C \cdot [(b-x)(x-a)]^{1/2},$ as functions of $x\in (a,b),$ for $a=-1, \  b=1,$
and $C= 1.35;$
the two curves appear to be almost undistinguishable.

\begin{figure}
\centering
\includegraphics[height=0.33 \textheight]{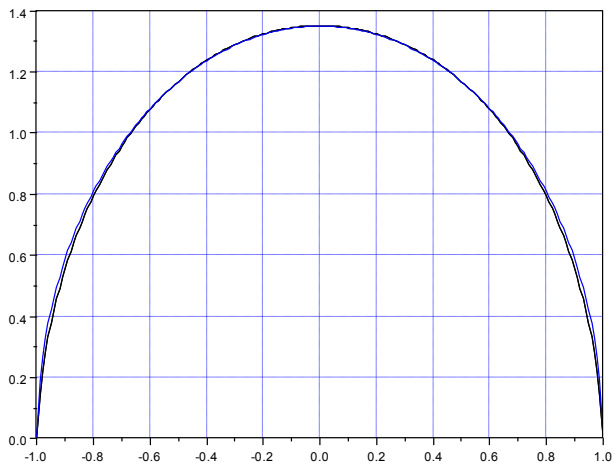}
\caption{Plots of the mean exit time, $E( \tau _{-1, 1} (x, 0)),$ of integrated BM from the interval $(-1,1)$
(lower curve), and of the function $z(x) = 1.35 \cdot (1-x^2)^{1/2}$ (upper curve), as functions of $x \in (-1,1).$
}
\end{figure}
We have also calculated the second order moment of the first-exit time of integrated BM,
by summing the first $20$ addends of the
series in \eqref{explictauquadIBM}.
In the Figure 2, we plot $E \left [ \left (\tau _{a, b} (x, 0) \right ) ^2 \right ], \ E ^2\left [ \tau _{a, b} (x, 0)  \right ]$
and the variance \par\noindent
$Var \left [ \tau _{a, b} (x, 0) \right ]=
E\left [ \left (\tau _{a, b} (x, 0) \right ) ^2 \right ] - \left ( E  \left [ \tau _{a, b} (x, 0) \right ] \right ) ^2 ,$ as a function of $x \in (-1, 1),$
for $a =-1, \ b =1;$ as we see, the maximum  of $Var \left [ \tau _{a, b } (x, 0) \right ]$ is about $10 \% $ times
the maximum  of $E( \tau _{-1, 1} (x, 0)).$
\par\noindent

\begin{figure}
\centering
\includegraphics[height=0.33 \textheight]{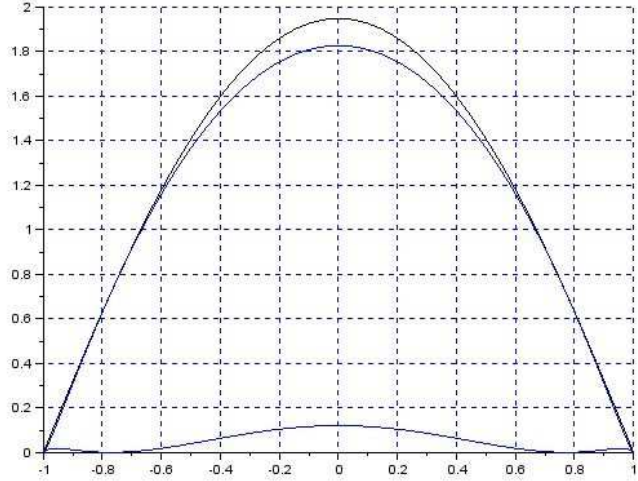}
\caption{From top to bottom: plot of the second moment (first curve), the square of the first moment (second curve), and the variance of the first-exit time $ \tau _{-1, 1} (x, 0)$  (third curve) of integrated BM from the interval $(-1,1),$
as functions of $x \in (-1.1).$
}
\end{figure}

\bigskip

As in the one boundary case, if we introduce a randomness in the
starting point, replacing $X(0)=x \in (a,b)$ with a random
variable $\eta ,$ having  density $g(x)$ whose support is the
interval $(a,b),$ we can consider the average  exit time
over all initial positions $\eta \in (a,b).$ If
$y=0,$ this quantity is:
$$ \overline T _ {a,b}= \int _ a ^b E( \tau _{a,b} (x,0)) g(x) dx .$$
In the case of integrated BM, $\overline T _ {a,b}$ can be
calculated by using the expression of $E( \tau _{a,b} (x,0))$
given by \eqref{explictauIBM}. We obtain:
\begin{equation} \label{TBM}
\overline T _ {a,b} = \frac {3^ {1/3} 2^{7/3} \Gamma (\frac 4 3 ) (b-a) ^ {2/3}  } {\pi ^ {5/3} }
\sum _ {k=0} ^ \infty (-1)^k \frac 1 {(2k+1) ^ {5/3} } \int _a ^b \cos \left [ \frac { \pi (2k+1)} {b-a}  \left (x - \frac {a+b} 2 \right )  \right ] g(x)dx
\end{equation}
(it has been possible to exchange the integral of the sum with the sum of the integrals, thanks to the dominated convergence theorem);
the integral in \eqref{TBM} equals $E(U_k),$ where
$ U_k=  \cos \left [ \frac {\pi (2k+1) } { b-a} \left ( \eta -
\frac { a+b} 2 \right ) \right ] \le 1 .$  Therefore:

\begin{equation} \label{averagemeanexittimegeneralg}
\overline T _ {a,b} = \frac {3^ {1/3} 2^{7/3} \Gamma (\frac 4 3 ) (b-a) ^ {2/3}  } {\pi ^ {5/3} }
\sum _ {k=0} ^ \infty (-1)^k \frac 1 {(2k+1) ^ {5/3} } E(U_k).
\end{equation}

In the special case when $g$ is the uniform density in
the interval $(a,b),$   we get by calculation:
$$  \overline T _ {a,b} =
\frac {3^ {1/3} 2^{7/3} \Gamma (\frac 4 3 ) (b-a) ^ {2/3}  } {\pi ^ {5/3} }
\sum _ {k=0} ^ \infty (-1)^k \frac 1 {(2k+1) ^ {5/3} } \int _a ^b \cos \left [ \frac { \pi (2k+1)} {b-a}  \left (x - \frac {a+b} 2 \right )  \right ]
\ \frac 1 { b-a} \ dx $$
\begin{equation} \label{averagemeanexittime}
= \frac {3^ {1/3} 2^{10/3} \Gamma (\frac 4 3 ) (b-a) ^ {2/3}  } {\pi ^ {8/3} } \sum _{k=0} ^ \infty \frac 1 { (2k+1)^ {8/3} } \ .
\end{equation}
Thus, $\overline T _ {a,b}= const \cdot (b-a)^ {2/3}.$ This confirms the result by
Masoliver and Porr\`a (see \cite{masoliver:phyrev96},
\cite{masoliver:phyrev95}), obtained by the Kolmogorov's equations
approach in the case of integrated BM, with $y=0$ and uniform distribution of the $X-$
starting point, according to which, the dependence of $\overline T _ {a,b}$ on the size $L=(b-a)$ of the interval, is $L^{2/3}.$
\bigskip

As far as integrated OU process is concerned, the moments of $ \tau _ {a,b} (x,0)$  can be found again by  formula \eqref{tauIBM},
where $\widehat \rho (t)$ can be deduced from
Example 2; however,
it is not possible to calculate explicitly the integral which appears in the expression of $A_k(x),$ so it has to be numerically computed. Since the integrand function decreases
exponentially fast, it suffices to calculate the integral over the interval $(0, 10),$ to obtain precise enough estimates. In the Figure 3
we have plotted, for comparison, the numerical evaluation
of the mean exit time of integrated OU process with $y= \beta =0 , $ from the interval $(-1,1),$
as a function of $x \in (-1,1),$ for $\sigma =1 $ and several values of $\mu ;$ in the Figure 4
we we have plotted the numerical evaluation
of $E \left [ \left (\tau _{- 1, 1} (x, 0) \right ) ^2 \right ], \ E ^2\left [ \tau _{- 1, 1} (x, 0)  \right ]$
and the variance
$Var \left [ \tau _{- 1, 1} (x, 0) \right ]$ of the first exit time of integrated OU process, for $\sigma =1 $ and $\mu=1.$ As we see,
the maximum of
$Var \left [ \tau _{- 1, 1} (x, 0) \right ] $ is about
$5 \% $ times the maximum of $E  \left (\tau _{- 1, 1} (x, 0) \right ).$

\begin{figure}
\centering
\includegraphics[height=0.37 \textheight]{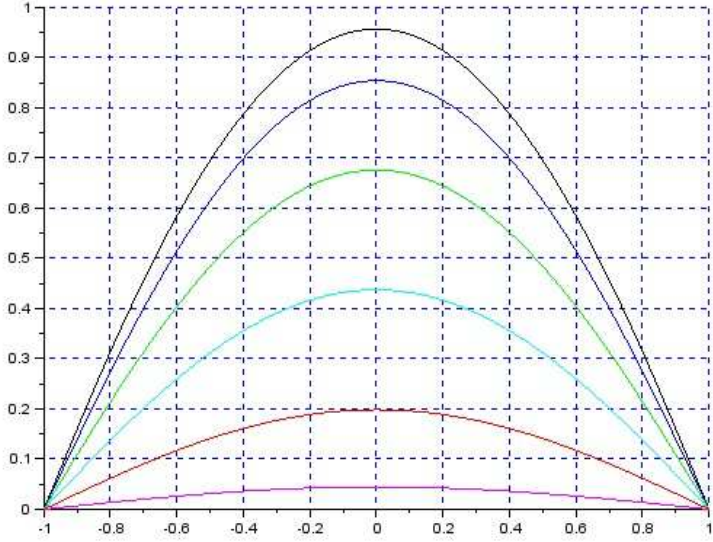}
\caption{Plot of numerical evaluation of the mean exit time,
$ E \left ( \tau _{-1, 1} (x, 0) \right ),$ of integrated OU
with $\beta = y =0, $ from the interval $(-1,1),$ as a function
of $x \in (-1,1),$ for $ \sigma =1$ and several values of  $\mu.$
From top to bottom, with respect to the peak of the curve: $\mu = 2; 1.8; 1.6; 1.4; 1.2; 1.$
}
\end{figure}

\begin{figure}
\centering
\includegraphics[height=0.37 \textheight]{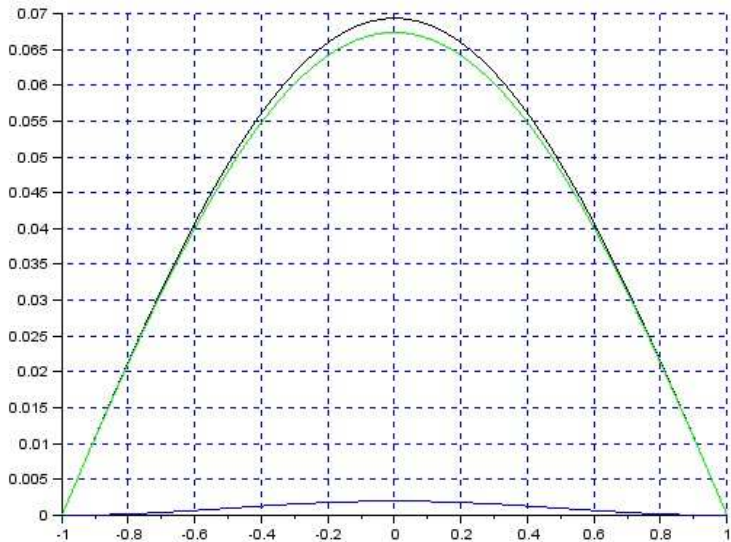}
\caption{
From top to bottom: plot of the second moment (first curve), the square of the first moment (second curve),
and the variance of the first-exit time $ \tau _{-1, 1} (x, 0)$  (third curve) of integrated OU
with $y= \beta =0 , $ from the interval $(-1,1),$ as a function of $x \in (-1,1),$ for $ \sigma = 1, \mu =1.$
}
\end{figure}

\bigskip

Finally,  we mention  the exit probabilities of the integrated Gauss-Markov process $X$ through  the ends of the interval $(a,b),$
namely:
$$\pi _a (x,y) = P \left ( \tau _a (x,y) < \tau _ b (x,y) \right )= P \left (X( \tau _ {a,b} (x,y)) =a \right ),$$ and
$$\pi _b (x,y) = P \left ( \tau _b (x,y) < \tau _ a (x,y) \right ) = P \left (X( \tau _ {a,b} (x,y)) =b \right ).$$
Recalling the well-known formulae
for exit probabilities of BM, we get, for $y=0$ and $x \in (a,b):$
$$ \pi _a (x,0)= P \left ( x+ \widehat B ( \widehat \tau _{a,b} (x,0)) =a \right )= \frac { b-x} {b-a } , \ \
\pi _b (x,0)= P \left ( x + \widehat B ( \widehat \tau _{a,b} (x,0)) =b  \right )= \frac { x-a} {b-a }.
$$
Notice that, in the case of integrated BM, several probability laws related to the couple $ \left ( \tau _{a,b}, \ B _ {\tau _{a,b}} \right )$ were evaluated in
\cite{lachal:crasp97} (in particular,
explicit formulae for $\pi _a (x,0)$ and $\pi _b (x,0)$ were obtained),
but they are written in terms of special functions.



\begin{thebibliography}{spc}

\bibitem [1] {abundo:smj13}
Abundo, M., 2013. \newblock
On the representation of an integrated Gauss-Markov process.
\newblock{Scientiae Mathematicae Japonicae Online} e-2013, 719–-723.

\bibitem [2] {abundo:stapro12}
Abundo, M., 2012. \newblock
An inverse first-passage problem for one-dimensional diffusions with random
      starting point.
\newblock{Statistics  and Probability Letters} 82 (1), 7--14.

\bibitem [3] {abundo:ijam12}
Abundo, M., 2012. \newblock
First-passage time of a stochastic integral process trough a linear boundary.
\newblock{International Journal of Applied  Mathematics (IJAM)} 25 (1), 41--49.


\bibitem [4] {abundo:ija08}
Abundo, M., 2008. \newblock
On the distribution of the time average of a  jump-diffusion process.
\newblock{International Journal of Applied  Mathematics (IJAM)} 21 (3), 447--454.

\bibitem [5] {abundo:stapro02}
Abundo, M., 2002. \newblock
Some conditional crossing results of Brownian motion over a piecewise-linear
boundary.
\newblock{Statistics and Probability Letters} 58 (2), 131--145.

\bibitem [6] {abundo:ric01}
Abundo, M., 2001. \newblock
Some results about boundary crossing for Brownian motion.
\newblock{Ricerche di
 Matematica} L (2), 283--301.

\bibitem [7] {ali:10}
Alili, L., and Patie, P., 2010. \newblock
Boundary-crossing identies for diffusions having the time-inversion property.
\newblock{J. Theor. Probab.} 23, 65--84.

\bibitem [8] {and:01}
Andersen, T.G., Bollerslev, T., Diebold, F.X. and Labys, P., 2001.
\newblock The distribution of realized exchange rate volatility.
\newblock{J. Amer. Statist. Assoc.} 96, 42--55.

\bibitem [9] {bache:41}
Bachelier, L., 1941. \newblock
Probabilit\'es des oscillation maxima.
\newblock{C. R. Acad. Sci. Paris} 212, 836--838.

\bibitem [10] {ben:13}
Benedetto, E., Sacerdote, L. and
Zucca, C., 2013.
\newblock A first passage problem for a bivariate diffusion process: Numerical solution with an application to neuroscience when the process is Gauss–Markov.
\newblock{J. Comput. Appl. Math.} 242, 41--52.

\bibitem [11] {bre:67}
Breiman, L., 1965/1966. \newblock
First exit times from a square root boundary.
In: Proc. Fifth Berkeley Sympos. Math. Statis. and Probability, Berkeley, Calif. Contributions to Probability
Theory. Part 2, vol II, 9--16. \newblock
University of California Press, Berkeley.


\bibitem [12] {darling:ams53}
Darling, D. A. and Siegert, A.J.F., 1953. \newblock
The first passage problem for a continuous Markov process.
\newblock   {Ann. Math. Statis.}, 24, 624--639.

\bibitem [13] {dan:ast82}
Daniels, H. E., 1982. \newblock
Sequential tests constructed from images.
\newblock   {Ann. Statist.}, 10 (2), 394--400.

\bibitem [14] {dan:jap69}
Daniels, H. E., 1969. \newblock
Minimum of a stationay Markov process superimposed on a U-shaped trend.
\newblock   {J. Appl. Probab.}, 6 (2), 399--408.

\bibitem[15] {doob:49}
Doob, J.L., 1949. \newblock
Heuristic approach to the Kolmogorov-Smirnov theorem. \newblock
{Ann. Math. Statist.}, 20, 393--403.

\bibitem [16] {gen:00}
Genon-Catalot, V., Jeantheau, T. and Lar\'edo, C., 2000.
\newblock Stochastic volatility models as hidden Markov models and
statistical applications.
\newblock{Bernoulli} 6, 1051--1079.

\bibitem [17] {gold:ams71}
Goldman, M., 1971. \newblock
On the first-passage time of the integrated Wiener process.
\newblock{Ann. Math. Statis.} 42 (6), 2150--2155.

\bibitem [18] {gork:ams75}
Gor'kov, Ju P., 1975. \newblock
A formula for the solution of a certain boundary value problem for the stationary equation of
Brownian motion.
\newblock{Probab. Theory. Relat. Fields} 81 (1), 79--109.

\bibitem [19] {glot:00}
Gloter, A., 2000. \newblock Parameter estimation for a discrete
sampling of an integrated Ornstein-Uhlenbeck process.
\newblock{Statistics} 35, 225--243.

\bibitem [20] {groe:ptr89}
Groeneboom, P., 1989. \newblock
Brownian motion with a parabolic drift and Airy functions.
\newblock{Sov. Math. Dokl.} 16, 904--908.

\bibitem [21] {hesse:iamsa05}
Hesse, C.H., 2005. \newblock
On the first-passage time of Integrated Brownian motion.
\newblock{Journal of Applied  Mathematics and Stochastic Analysis} 3, 237--246.

\bibitem[22]
{kle:stoch05}
Klebaner, F.C., 2005. \newblock
Introduction to Stochastic Calculus with Applications. Second Edition. \newblock
Imperial College Press, London.

\bibitem[23]
{kolmogorov:anm34}
Kolmogorov, A., 1934 \newblock
Zuf$\ddot{a}$llige Bewegungen (zur Theorie der Brownschen Bewegung) \newblock
\newblock{Ann. of Math.} (2) 35, no. 1, 116-117 (German).

\bibitem [24] {lachal:crasp97}
Lachal, A., 1997. \newblock
Temps de sortie d'un intervalle born\'e pour
l'int\'egrale du mouvement Brownien.
\newblock{C.R. Acad. Sci. Paris} 324, serie I, 559--564.

\bibitem [25] {lachal:jap93}
Lachal, A., 1993. \newblock
L'integrale du mouvement Brownien.
\newblock{J. Appl. Prob.} 30, 17--27.

\bibitem [26] {lachal:aihp91}
Lachal, A., 1991. \newblock
Sur le premier instant de passage de l'integrale du mouvement Brownien.
\newblock{Annales de l' I.H.P.} B, 27 (3), 385--405.

\bibitem [27] {lansky:89}
Lansky, P. and Smith, C.E., 1989.\newblock
The effect of a random initial value in neural first-
passage-time models.
\newblock{Math. Biosci.} 93 (2), 191–-215.

\bibitem [28] {lefebvre:sjam89}
Lefebvre, M., 1989. \newblock
First-passage densities of a two-dimensional process.
\newblock{SIAM J. Appl. Math.} 49 (5), 1514--1523.

\bibitem [29] {lefebvre:spa89}
Lefebvre, M., 1989. \newblock
Moment generating function of a first hitting place for the integrated Ornstein-Uhlenbeck process.
\newblock{Stoch. Proc. Appl.} 32, 281--287.

\bibitem[30]
{marcus:mar06}
Marcus M.B. and Rosen, J., 2006. \newblock
Markov processes, Gaussian processes, and local time. \newblock
Cambridge University Press, Cambridge.

\bibitem [31] {masoliver:phyrev96}
Masoliver, J. and Porr\`a, J.M., 1996. \newblock
Exact solution to the exit-time problem for an undamped free particle driven by Gaussian white noise.
\newblock{Phys. Rev. E} 53 (3), 2243--2256.

\bibitem [32] {masoliver:phyrev95}
Masoliver, J. and Porr\`a, J.M., 1995. \newblock
Exact solution to the mean exit-time problem for free inertial processes driven by Gaussian white noise.
\newblock{Phys. Rev. Lett.} 75 (2), 189--192.

\bibitem [33] {mehr:65}
Mehr, C.B. and McFadden, J.A., 1965. \newblock
Certain properties of Gaussian processes and their first-passage times.
\newblock{J. R. Statist. Soc. B} 27, 505--522.

\bibitem [34] {mckean:jmku63}
McKean, H.P. Jr., 1963. \newblock
A winding problem for a resonator driven by a white noise.
\newblock{J. Math. Kyoto Univ.} 2 (2), 227--235.

\bibitem[35]
{ric:smj08}
Nobile, A.G., Pirozzi, E., Ricciardi, L.M. , 2008. \newblock
Asymptotics and evaluations of FPT densities through varying boundaries for Gauss-Markov processes. \newblock
{Scientiae Mathematicae Japonicae} 67, (2), 241--266.

\bibitem[36]
{revuzyor:con91}
Revuz, D. and Yor, M., 1991. \newblock
Continous martingales and Brownian motion. \newblock
Springer-Verlag, Berlin Heidelberg.

\bibitem[37]
{ross:pmodel10}
Ross, S.M., 2010. \newblock
Introduction to Probability Models. Tenth Edition. \newblock
Academic Press, Elsevier, Burlington.

\bibitem [38] {salm:aap88}
Salminen, P., 1988. \newblock
On the first hitting time and the last exit time for a
Brownian motion to/from a moving boundary.
\newblock{Adv. Appl. Probab.} 20 (2), 411--426.

\bibitem [39] {sato:77}
Sato, S., 1977. \newblock
Evaluation of the first-passage time probability to a square root boundary for the Wiener process..
\newblock{J. Appl. Probab.} 14 (4), 850--856.







\end{thebibliography}
\end{document}